\documentclass{article}
\newcommand{\To}{\Rightarrow}
\newcommand{\ga}{\Gamma}
\newcommand{\de}{\Delta}
\newcommand{\skep}{\hskip2em}
\usepackage{irule}
\usepackage{amssymb,verbatim}
\usepackage{array,amsmath}
\newtheorem{theorem}{Theorem}[section]
\newtheorem{definition}{Definition}[section]
\newtheorem{lemma}{Lemma}[section]
\newtheorem{corollary}{Corollary}[section]

\newcommand{\qed}{\hfill $\Box$}

\newcommand{\G}{Grzegorczyk}
\newcommand{\cd}{{\bf CD}}
\newcommand{\LD}{$L(D)$}
\newcommand{\forces}{\Vdash_{\cal M}}
\newcommand{\forcesi}{\Vdash_i}
\newcommand{\forcesj}{\Vdash_j}

\begin{document}

\title{Failure of interpolation in the intuitionistic logic of constant domains}

\author{Grigori Mints\\
Stanford University \and
    Grigory Olkhovikov\\
Ural Federal University \and Alasdair Urquhart\thanks{\ Research
supported by the Natural Sciences and Engineering Research Council of Canada.}\\
University of Toronto}

\maketitle

\begin{abstract}

This paper shows that the interpolation theorem fails in the
intuitionistic logic of constant domains. This result refutes two
previously published claims that the interpolation property holds.

\end{abstract}

\section{Introduction}

The main result of this paper demonstrates that the interpolation theorem fails in the
intuitionistic logic of constant domains.
This shows that two previously published proofs of
the interpolation property are incorrect.

The model theory for the intuitionistic logic of constant domains
(\cd) was proposed by \G\ in a paper \cite{Grzegorczyk1964} of
1964, inspired by Paul Cohen's notion of forcing, as well as the
earlier semantics for intuitionistic logic proposed by Evert
W.~Beth \cite{Beth1956,Beth1959}. The logic resulting from the
semantics is stronger than intuitionistic predicate logic, since
it includes the scheme
\[
(D) \:\:\:\: \forall x (A \vee B) \rightarrow (A \vee \forall x
B),
\]
where $x$ is not free in $A$; this scheme is not provable
intuitionistically. \G\ proposed the semantics as a
``philosophically plausible formal interpretation of
intuitionistic logic,'' independently of the near-contemporary
work of Saul Kripke \cite{Kripke1965}. \G\ observes that his
semantics validates the schema (D), but goes on to propose a
modification to the forcing relation for disjunctions and
existential formulas that gives an exact interpretation for
intuitionistic predicate logic.

Sabine G\"{o}rnemann \cite{Gornemann1971} proved  that the
addition of the scheme to intuitionistic predicate logic is
sufficient to axiomatize \G's logic; D.~Klemke \cite{Klemke1971}
and Dov Gabbay \cite{Gabbay1969} gave independent proofs of the
same result.

There are two published proofs that the interpolation theorem
holds for the intuitionistic logic of constant domains. The first
is by Dov M.~Gabbay \cite{Gabbay1977}, who gives a model-theoretic
proof. The second is by L\'{o}pez-Escobar \cite{LopezEscobar1981},
using proof theory. In a later article  \cite{LopezEscobar1983},
L\'{o}pez-Escobar admitted that there was an error in his earlier
proof of interpolation.

The first author of the present paper wrote in his review of
\cite{LopezEscobar1983} in the Russian review journal `Mathematika'
[RZhmat], 7A70, July 1984: ``The first proof of the interpolation
theorem for \cd\ given by Dov Gabbay (reviewed in RZhmat,
11A83, November 1971) contained gaps and specialists have
different opinions on whether these gaps were closed later. The
proof from the previous paper by the present author [L.-E.]
(reviewed in RZhmat, 12A40, December  1981) consisted in essence
of the reference to a cut-elimination theorem for some
Gentzen-type axiomatization for \cd. In the paper under review its
author reproduces a counterexample to this claim given by M.
Fitting. (Another counterexample was given by the present reviewer
in RZhmat, 12A40, December  1981.) Therefore the question of the
validity of the interpolation theorem for \cd\ is not clear. The present
author [L.E.] tries to prove that \cd\ does not have any cut-free
Gentzen-style axiomatization with a finite number of rules where
formulas are analyzed only up to a finite depth $k$. He concludes
this from the statement of necessity in such axiomatics of a derivable formula
 \[
 \forall x(p\vee B(x))\to (p\vee (q\to\forall x_1\ldots x_{k-1}\forall x B(x)))
 \]
However this formula is derivable in the axiomatics obtained by adding
the rule
 \[
 \infer{\forall x(p\vee Bx),\Gamma\To C}  {p,\Gamma\To C &\forall xBx,\Gamma\To C}
\]
to ordinary intuitionistic rules. The error by the author is in
the case C2 which is not treated.''

In his paper of 1977, Gabbay shows that the strong Robinson
consistency theorem does not hold for \cd; however, he goes
on to claim that the weak Robinson consistency theorem holds for
\cd, and from this deduces the interpolation theorem. The
first (negative) result appears in the comprehensive monograph
co-authored by Gabbay and Maksimova \cite{GabbayMaksimova}, but
the proof of interpolation does not appear.

Our non-interpolation result can be used to justify the impression of
E. G. K. L\'{o}pez-Escobar in \cite{LopezEscobar1981}  that {\bf
CD} does not admit a ``good'' cut-free formulation. Recall that
one of the familiar proofs of interpolation for classical
predicate logic  (say in \cite{Takeuti1987}) goes by induction on
a cut-free Gentzen-style derivation: interpolants are defined for
axioms, and simple prescriptions are given  for transferring
interpolants from premises of an  inference to its conclusion.
Hence our result shows that \cd\ does not admit a formulation
where interpolants are defined for axioms and there is an explicit
definition of the interpolant for the conclusion of every
inference  rule from interpolants for the premises of the rule.

In fact this idea motivated  the method by which the first author of the
present paper constructed the counterexample to interpolation
treated below. The implication $\ga\to \de$ lacking an interpolant
was designed in such a way that any obvious (to the author) proof
of it would involve a cut over a formula connecting predicates to
be eliminated with the schema (D) above.

An important non-interpolation result was obtained  in
\cite{Fine1979}. Kit Fine proved there that any normal predicate
modal logic with Barcan formula as an axiom between {\bf K} and
{\bf S5} lacks interpolation. One can try to get our result from
this using the G\"{o}del-Tarski translation $^\Box$ of
intuitionistic logic into {\bf S4} which just prefixes  the
necessity symbol $\Box$ to every subformula. The result $F^\Box$
of this translation is S4-valid (derivable) if and only if $F$ is
valid (derivable) intuitionistically. Moreover the translation is
sound and faithful also for constant domains, that is as an
operation from \cd\ into the modal predicate logic based on {\bf
S4} plus the Barcan formula. There are  however essential
obstacles blocking an attempt to reduce the problem for \cd\ to
similar problem for modal systems.  First, counterexamples to
interpolation used by Kit Fine do not translate in an obvious way
back from modal logic to intuitionistic logic or \cd. Second, our
counterexample $\ga\to\de$  {\em does} in fact have a modal
interpolant. More about this in the section \ref{scounterexample}.
Third, Fine refutes not only Craig interpolation, but also Beth
definability theorem for his modal systems. It is not known at
this moment whether \cd\ has the Beth property, although this does
not look plausible.

A review of \cite{Fine1979} by Saul Kripke \cite{Kripke1983}
analyzed Fine's construction in terms of second order
quantification of ``redundant'' variables and formed a background
for computations in our section \ref{scounterexample}.

The counterexample $\ga\to\de$ to interpolation for \cd\ given
at the beginning of that section  was constructed  by the first
author at the beginning of 2008. The third author joined the
efforts to prove that no interpolant exists for the example in 2011 and
computed ``first order equivalents'' of $\exists R\ga$ and
$\forall S\de$ used in our proof. (These results can be stated in
a language recalling the language of hybrid modal logic
\cite{BlackburndeRijkeVenema2001}.) After that our efforts were
stalled till the second author arrived in Stanford and joined the
project. This resulted in a construction of models proving an
interpolant for $\ga\To\de $ is indeed impossible based on his
characterization of first-order  translations of intuitionistic
formulas \cite{Olkhovikov2012}.

In Section \ref{smodeltheory} we recall the definition of \cd\
and its model theory.
In Section \ref{scounterexample} we present our counterexample and
compute ``candidate interpolants'' in the first order language of
Kripke-style models.
Section \ref{sasimulations} contains the definition of asimulation
(a generalization of bisimulation for modal logic) due to the
second author and proof of the easy part (necessity) of
asimulation condition needed for our result. Section
\ref{srefuting} contains the construction of the pair of models
proving non-interpolation, and a detailed proof of the failure
of interpolation for \cd.

\section{Model theory}\label{smodeltheory}

We formulate the language $L$ of the logic \cd\ using the
propositional connectives $\wedge$, $\vee$, $\rightarrow$ and
$\bot$, together with the universal and existential quantifiers
$\forall$ and $\exists$; the negation operator $\neg A$ is defined
as $A \rightarrow \bot$. The atomic formulas are of the form
$Px_1, \dots, x_k$; we allow the case where the sequence of
variables $x_1, \dots, x_k$ is empty, that is to say, where $P$ is
a propositional variable. We use the notation $A[x_1, \dots,
x_k]$, or $A[\vec{x}]$, for a formula of the language $L$, where
all the free variables in the formula appear in the sequence $x_1,
\dots, x_k$ (some of the variables in the sequence may not appear
in the formula). We employ the notation $L(P,Q)$ for the
sublanguage of $L$ in which the only predicate symbols are the
one-place predicates $P$ and $Q$.

If $D$ is a non-empty set, then the language \LD\ is
obtained from $L$ by adding a distinct constant $\mathbf a$ for
every element $a$ in $D$. If $A[{\mathbf a_1}/x_1, \dots, {\mathbf
a_k}/x_k]$ is a sentence in the expanded language \LD, where
$A[x_1, \dots, x_k]$ is a formula of $L$, then we define
$A[{\mathbf a_1}/x_1, \dots, {\mathbf a_k}/x_k]$ to be a theorem
of \cd\ if the universal closure $\forall x_1 \dots \forall x_k
A[x_1, \dots x_k]$ is a theorem of \cd.

\G's model theory is very similar to that of Kripke
\cite{Kripke1965}, but simpler. Kripke's model theory involves a
quasi-ordered set of classical models, where the domains of the
models can expand along the quasi-ordering. By contrast, \G's
models have a fixed, constant domain. Thus the \G\ model theory
represents a static ontology, whereas Kripke's represents an
expanding ontology, where new objects can be created as knowledge
grows.

A \G-model, or {\em G-model} for short, is a structure
\[
{\cal M} = \langle W, \leq, \mathbf{w}, D, \phi \rangle,
\]
 where $W$ is a non-empty set (the set of information states in the terminology of \G), $\leq$
is a quasi-ordering on $W$ (a reflexive transitive relation on
$W$), $\mathbf{w}$ is the base state in the model satisfying
$\forall v  \in W ( \mathbf{w} \leq v)$, $D$ is a non-empty set,
and for each $k$-ary predicate symbol $P$ in $L$, $\phi(P)$ is a
$k+1$-ary relation contained in $W \times D^k$ that satisfies the
monotonicity condition
\[
 [ v \leq w \wedge \langle v, a_1, \dots, a_k \rangle \in \phi(P) ] \Rightarrow \langle w, a_1, \dots, a_k \rangle \in \phi(P).
\]
The forcing relation $\forces$ holds between states in $W$ and
sentences in \LD\ by the following inductive definition:
\begin{enumerate}
\item
  $v \forces P{\mathbf a}_1, \dots, {\mathbf a}_k \Leftrightarrow \langle v, a_1, \dots, a_k \rangle \in \phi(P)$,
\item
  $v \forces A \wedge B \Leftrightarrow (v \forces A \wedge v \forces B)$,
\item
  $v \forces A \vee B \Leftrightarrow (v \forces A \vee v \forces B)$,
\item
 $v \forces A \to B \Leftrightarrow \forall w \geq v (w \forces A \to w \forces B)$,
\item
  $v \forces \bot$ never holds,
\item
  $v \forces \exists x A \Leftrightarrow \exists a \in D ( v \forces A[{\mathbf a}/x]) $,
\item
  $v \forces \forall x A \Leftrightarrow \forall a \in D ( v \forces A[{\mathbf a}/x]) $.
\end{enumerate}

It is easy to verify that this semantics is sound for \cd, in the
sense that if $v$ is a state in a \G-model, and $A$ a sentence of
\LD\ that is a theorem of \cd, then $v \forces A$. The
completeness theorem for \cd\ asserts that if $A$ is a sentence
that is not a theorem of \cd, then there is a \G-model $\cal M$ with base
point $\mathbf{v}$ so that $\mathbf{v} \not\forces A$.

The semantics is easily extended to a second-order version, where
the second-order variables of arity $k$ range over $k+1$-ary
relations $R$ over $W \times D^k$ that satisfy a version of the
monotonicity condition above:
\[
 [ v \leq w \wedge \langle v, a_1, \dots, a_k \rangle \in R ] \Rightarrow \langle w, a_1, \dots, a_k \rangle \in R.
\]
We shall use the second-order semantics only in the one-quantifier form.

\section{The Counterexample}\label{scounterexample}

In this section, we produce the counterexample used in refuting
the interpolation theorem for \cd. The two formulas forming the
antecedent and consequent of the implication for which no
interpolant exists in \cd\ are as follows:
\[
 \Gamma = [\forall x \exists y (Py \wedge (Qy \rightarrow Rx)) \wedge \neg \forall x Rx],
\]
\[
 \Delta = \forall x (Px \rightarrow (Qx \vee S)) \rightarrow S.
\]

\begin{lemma}\label{Implicationlemma}
The implication $\Gamma \rightarrow \Delta$ is valid in all
\G-models.
\end{lemma}

{\bf Proof.} We give two proofs, the first model-theoretic,
the second a deductive proof in a sequent calculus.

First, we give a model-theoretic argument for the Lemma.
Let $v$ be a state in a \G-model ${\cal M}$ with domain $D$ so
that $v \forces \Gamma$; we wish to show that $v \forces \Delta$.
Let $w$ be a state in $\cal M$ so that $v \leq w$, where $w
\forces  \forall x (Px \rightarrow (Qx \vee S))$. For an arbitrary
$a \in D$, since $v \forces \Gamma$, $w \forces \exists y ( Py
\wedge (Qy \rightarrow Ra ))$, so that for some $b \in D$, $w
\forces Pb \wedge (Qb \rightarrow Ra)$. We also have $w \forces Pb
\rightarrow ( Qb \vee S)$, hence $w \forces (Qb \vee S)$ and so $w
\forces Ra \vee S$. Since $a$ was an arbitrary element of $D$, it
follows that $w \forces \forall x (Rx \vee S)$. Hence, by $(D)$,
$w \forces \forall x Rx \vee S$. Since $v \forces \Gamma$, it
follows that $w \forces \neg \forall x Rx$, so $w \forces S$,
showing that $v \forces \Delta$.

It follows by the completeness theorem for \cd\ that $\Gamma
\rightarrow \Delta$ is a theorem of \cd. However, we can argue
directly for this in a multiple succedent sequent calculus
formulation of \cd  (see, e.\,g. \cite{Mints2000}), using the
notation $A_1\ldots A_n \To B_1\ldots B_m$ for sequents.

Employing the abbreviation
\[
\alpha:\equiv\ \forall x (Px \rightarrow (Qx \vee S)),
\]
we derive intuitionistically a nucleus of our implication.
\[
\infer{\forall x \exists y (Py \wedge (Qy \rightarrow Rx)),\alpha
\To\forall x(Rx\vee S)}
  {\infer{\forall x \exists y (Py \wedge (Qy \rightarrow Rx)),\alpha \To Rx\vee S}
     {\infer{ \exists y (Py \wedge (Qy \rightarrow Rx)),\alpha \To Rx\vee S}
        {\infer{ (Py \wedge (Qy \rightarrow Rx)), \forall x (Px \rightarrow (Qx \vee S)) \To Rx\vee S}
            {\infer{ (Py \wedge (Qy \rightarrow Rx)),(Py \rightarrow (Qy \vee S)) \To Rx\vee S}
               {\infer{ Py, (Qy \rightarrow Rx), (Py \rightarrow (Qy \vee S)) \To Rx\vee S}
                 {\infer{ (Qy \rightarrow Rx), (Qy \vee S) \To Rx\vee S\skep  Py\To Py}
                  { (Qy \rightarrow Rx),Qy  \To Rx\vee S\skep S\To Rx\vee S}
                 }
               }
            }
        }
   }
  }
\]

After that an application of
the cut rule with an instance of the axiom $D$ produces
the required result. We use the abbreviation

\[
\beta:= \forall x \exists y (Py \wedge (Qy \rightarrow Rx)).
\]

\[
\infer{\ga\To\de}
 {\infer{\ga,\alpha\To S}
   {\infer{\beta, \neg\forall xRx,\alpha\To S}
    {\infer{\beta,\alpha \To \forall xRx, S}
      {\deduce{\beta,\alpha \To \forall x(Rx\vee S)}{\vdots}
       &
        \infer[D,Cut]{\forall x(Rx\vee S)\To\forall xRx,S}
          {\infer{\forall xRx\vee S\To\forall xRx, S}
            {\forall xRx\To\forall xRx & S\To S}
          }
     }
    }
   }
 }
\]

\qed

\begin{corollary}\label{Implicationcorollary}
 The second-order implication $\exists R \Gamma \rightarrow \forall S \Delta$ is valid in all G-models.
\end{corollary}

It is possible to shed some further light on Lemma
\ref{Implicationlemma} and Corollary \ref{Implicationcorollary}
through an analysis of the second-order sentences $\exists R
\Gamma$ and $\forall S \Delta$. In the following Lemma, we extract
the content of these two sentences in the form of semantical
conditions on G-models.

\begin{lemma}\label{Semanticalinterpolants}
Let $\cal M$ be a G-model with base point $\mathbf{v}$.
\begin{enumerate}
 \item $\mathbf{v} \forces \exists R \Gamma$ if and only if $\cal M$ satisfies the semantical condition
\[
 I(P,Q) \equiv \forall w \exists a ( \mathbf{v} \forces Pa \wedge w \not\forces Qa).
\]
\item $\mathbf{v} \forces \forall S \Delta$ if and only if $\cal
M$ satisfies the semantical condition
\[
 J(P,Q) \equiv \forall w \exists a ( w \forces Pa \wedge w \not\forces Qa).
\]
\end{enumerate}
\end{lemma}

{\bf Proof.} {\it Part 1} $(\Rightarrow)$. Assume that $\mathbf{v}
\forces \exists R \Gamma$; let $w$ be an arbitrary state in $\cal
M$. Then since $\mathbf{v} \forces \neg \forall x Rx$, $w
\not\forces \forall x Rx$, so that for some $b \in D$, $w
\not\forces Rb$. Since $\mathbf{v} \forces \Gamma$, $\mathbf{v}
\forces \exists y (Py \wedge (Qy \rightarrow Rb))$, so that for
some $a \in D$, $\mathbf{v} \forces Pa \wedge (Qa \rightarrow
Rb)$. It follows that $w \not\forces Qa$, showing that $\cal M$
satisfies $I(P,Q)$.

{\it Part 1} $(\Leftarrow)$. Conversely, assume that $\cal M$
satisfies $I(P,Q)$. Define $E \subseteq D$ as follows:
\[
 E := \{ a \in D \: | \:  \exists w \in {\cal M} ( \mathbf{v} \forces Pa \wedge w \not\forces Qa)\}.
\]
By assumption, $E$ is non-empty. Let $f$ be a surjective map from
$D$ to $E$, and define for any $a \in D$, $w \in W$,
\[
  w \forces Ra \Leftrightarrow w\forces Q(f(a)).
\]
Then for any $a \in D$, $\mathbf{v} \forces P(f(a))$, and
$\mathbf{v} \forces ( Q(f(a)) \leftrightarrow Ra)$, so that
$\mathbf{v} \forces \forall x \exists y (Py \wedge (Qy \rightarrow
Rx))$. Furthermore, for any $w \in {\cal M}$, there is a $b \in E$
so that $\mathbf{v} \forces Pb$ and $w \not\forces Qb$. Choose $a
\in D$ so that $f(a) = b$; we have $w \forces (Q(f(a))
\leftrightarrow Ra)$. Since $w \not\forces Q(f(a))$, $w
\not\forces Ra$. This shows that $w \not\forces \forall x Rx$, so
$\mathbf{v} \forces \neg \forall x Rx$

{\it Part 2} $(\Leftarrow)$. Assume that $\cal M$ satisfies
$J(P,Q)$. If $w$ is an arbitrary state in $\cal M$, where $w
\forces \forall x (Px \rightarrow (Qx \vee S))$, then by
assumption, there is an $a$ so that $w \forces Pa \wedge w
\not\forces Qa$. But $w \forces (Pa \rightarrow (Qa \vee S))$, so
$S$ must be true at $w$.

{\it Part 2} $(\Rightarrow)$. For the converse, let us assume that
$\cal M$ does not satisfy $J(P,Q)$. Then there must be a $w$ in
$\cal M$ so that
\[
 \forall a \in D ( w \forces Pa \Rightarrow w \forces Qa).
\]
Now define $S$ by setting $u \forces S$ if and only if $u > w$;
that is to say, the proposition $S$ holds exactly at those states
accessible from $w$, but not identical with $w$. We have to argue
that $w \forces \forall x (Px \rightarrow (Qx \vee S))$. Let $u$
be an arbitrary state accessible from $w$. If $u$ is $w$ itself,
we need to show that for any $a$, whenever $u \forces Pa$, then $u
\forces (Qa \vee S)$; this follows from our assumption. On the
other hand, if $u$ is not $w$, then $S$ holds at $u$, by
construction, so $(Qa \vee S)$ holds for any $a$. So, we've shown
that for any $a$, $w \forces (Pa \rightarrow (Qa \vee S))$, hence
$w \forces \forall x (Px \rightarrow (Qx \vee S))$. However, $S$
fails at $w$, by construction, so $\mathbf{v} \not\forces \forall
S \Delta$. \qed

\medskip

In notation of Lemma \ref{Semanticalinterpolants} note that
$J(P,Q)$ is true in all models of {\bf S4}+BF
(quantified {\bf S4} plus the Barcan formula)
if and only if
\[
S4+BF\models \Box\exists x(\Box Px\wedge \neg\Box Qx).
\]
Hence the implication $\Gamma\Rightarrow \Delta$ has a modal
interpolant. Theorem 5.1 below shows that this modal formula is
not equivalent to a G\"{o}del-Tarski translation (prefixing $\Box$
to all subformulas) of any intuitionistic formula.

Lemma \ref{Semanticalinterpolants} suggests a strategy to
demonstrate the failure of the interpolation theorem, namely by
showing that if $\Theta(\mathbf{v})$ is a first-order property of
G-models that is implied by $I(P,Q)$, and also implies $J(P,Q)$,
then it cannot be expressed by a formula of \cd\ involving only
the predicates $P$ and $Q$.

\section{Asimulations}\label{sasimulations}

In this section, we introduce the basic concept of {\em
CD-asimulation}, an asymmetric counterpart of the concept of
bisimulation familiar from the literature of modal logic
\cite[Chapter 2]{BlackburndeRijkeVenema2001}. The definition given
below can be considered as a more general version of the notion of
bisimulation for modal predicate logic defined by Johan van
Benthem \cite{VanBenthem2010}.

\begin{definition}\label{Asimulation}
If ${\cal M}_1 = \langle W_1, \leq_1, \mathbf{v}, D_1, \phi_1
\rangle$ and ${\cal M}_2 = \langle W_2, \leq_2, \mathbf{w}, D_2,
\phi_2 \rangle$ are G-models, then a {\em CD-asimulation} between
${\cal M}_1$ and ${\cal M}_2$ is a relation $Z$ satisfying the
following conditions:
\begin{enumerate}
 \item $Z \subseteq \bigcup_{k \geq 0} [ (W_1 \times D_1^k) \times (W_2 \times D_2^k)] \cup [ (W_2 \times D_2^k) \times (W_1 \times D_1^k)]$;
 \item $\{( v,\vec{d} Z w, \vec{e}) \wedge v \forcesi P[\vec{\bf d}] \} \Rightarrow w \forcesj P[\vec{\bf e}]$, for $P[\vec{x}]$ an atomic formula;
 \item $\{ (t,\vec{d} Z u,\vec{e}) \wedge u \leq_j v \} \Rightarrow
(\exists w \in W_i) ( t \leq_i w \wedge (w, \vec{d} Z v, \vec{e})
\wedge (v, \vec{e} Z w, \vec{d})  )$;
 \item $\{ t \in W_i \wedge (t,\vec{d} Z u, \vec{e}) \wedge f \in D_i \} \Rightarrow (\exists g \in D_j) ( t,\vec{d},f Z u, \vec{e},g )$;
 \item $\{ t \in W_i \wedge (t,\vec{d} Z u, \vec{e}) \wedge g \in D_j \} \Rightarrow (\exists f \in D_i) ( t,\vec{d},f Z u, \vec{e},g )$,
\end{enumerate}
where $\{i,j\} = \{1,2\}$, and $\forcesi$ and $\forcesj$ are the
forcing relations in ${\cal M}_i$ and ${\cal M}_j$.
\end{definition}

The concept of asimulation is due to the second author of this
paper. The version we are using here is a simplified version of
the general notion; the simplifications are possible because  we
are operating in the context of the constant domain semantics. The
second author has defined a more general version \cite{Olkhovikov2012}
that is suited to
the context of Kripke's semantics for intuitionistic predicate
logic, and has also proved that it serves to characterize the
first-order properties expressed by formulas of propositional and
predicate logic in that framework. In the present case, however,
we do not need the full characterization theorems; the following
Lemma is sufficient for our purpose of refuting the interpolation
theorem for the logic \cd.

\begin{lemma}\label{Soundness}
Let $Z$ be a CD-asimulation between ${\cal M}_1$ and ${\cal M}_2$,
and $A[\vec{x}]$ a formula of $L$. If $t,\vec{d} Z u, \vec{e}$ and
$t \forcesi A[\vec{\bf d}]$, then $u \forcesj A[\vec{\bf e}]$,
where $\{i,j\} = \{ 1,2 \}$.
\end{lemma}

{\bf Proof.} By induction on the complexity of the formula
$A[\vec{x}]$. For atomic formulas, the lemma holds by the
definition of CD-asimulation. The inductive steps for $\wedge$,
$\vee$ and $\bot$ are straightforward.

Now assume that $t,\vec{d} Z u, \vec{e}$ and that $t \forcesi
A[\vec{\bf d}] \rightarrow B[\vec{\bf d}]$. If $u \leq_j v$ and $v
\forcesj A[\vec{\bf e}]$, then by the third condition in
Definition \ref{Asimulation},
\[
 (\exists w \in W_i) ( t \leq_i w \wedge (w, \vec{d} Z v, \vec{e}) \wedge (v, \vec{e} Z w, \vec{d})  ).
\]
By inductive assumption, $w \forcesi A[\vec{\bf d}]$, so that $w
\forcesi B[\vec{\bf d}]$, since $t \forcesi A[\vec{\bf d}]
\rightarrow B[\vec{\bf d}]$. Again, by inductive assumption, $v
\forcesj B[\vec{\bf e}]$, showing that $u \forcesj A[\vec{\bf e}]
\rightarrow B[\vec{\bf e}]$.

Assume that $t,\vec{d} Z u, \vec{e}$ and that $t \forcesi \exists
x A[\vec{\bf d},x]$. Then $t \forcesi A[\vec{\bf d}, {\bf f}]$,
for some $f \in D_i$. By the fourth condition in Definition
\ref{Asimulation}, there is a $g$ in $D_j$ so that $t,\vec{d},f Z
u, \vec{e},g$. By inductive assumption, $u \forcesj A[\vec{\bf e},
{\bf g}]$, so that $u \forcesj \exists x A[\vec{\bf e}, x]$.

Assume that $t,\vec{d} Z u, \vec{e}$ and that $t \forcesi \forall
x A[\vec{\bf d},x]$. Let $g$ be an arbitrary individual in $D_j$.
By the fifth condition in Definition \ref{Asimulation}, $(\exists
f \in D_i) ( t,\vec{d},f Z u, \vec{e},g )$. Then $t \forcesi
A[\vec{\bf d}, {\bf f}]$, so by inductive assumption, $u \forcesj
A[\vec{\bf e}, {\bf g}]$, showing that $u \forcesj \forall x
A[\vec{\bf e}, x]$. \qed

\section{Refuting interpolation}\label{srefuting}

In this section, we define two G-models, ${\cal M}_1$ and ${\cal
M}_2$, together with a CD-asimulation between them (these
definitions and the ideas of the proofs given below are due to the
second author of this paper). This will enable us to carry out the
strategy outlined in the comments following Lemma
\ref{Semanticalinterpolants}. The states in both models are {\em
quasi-partitions}, given by the following definition. We use the
notation $\mathbb N$ for the set of positive natural numbers, and
for $k > 0$, $l \geq 0$, we write $k \mathbb{N} + l$ for the set
$\{ kn + l | n \in \mathbb{N} \}$, and $k \mathbb{N}$ for $k
\mathbb{N} + 0$.

\begin{definition}
\begin{enumerate}
\item A {\em quasi-partition} $(A,B,C)$ is defined by the
following conditions:
\begin{enumerate}
\item $A \cup B \cup C = \mathbb{N}$; \item $A,B,C$ are pairwise
disjoint; \item $A$ and $C$ are infinite; \item $B$ is either
empty or infinite.
\end{enumerate}
\item A quasi-order $\trianglelefteq$ on the set of all
quasi-partitions is defined by
\[
 (A,B,C) \trianglelefteq (D,E,F) \Leftrightarrow [ A \subseteq D \wedge F \subseteq C].
\]
\end{enumerate}
\end{definition}

We begin by defining two quasi-partitions that will serve as the
base points of our two models. The first is $\mathbf{v} =
(\mathbf{v}_1, \mathbf{v}_2, \mathbf{v}_3) = ( 3 \mathbb{N}, 3
\mathbb{N} + 1, 3 \mathbb{N} + 2 )$; the second is $\mathbf{w} =
(\mathbf{w}_1, \mathbf{w}_2, \mathbf{w}_3) = ( 2 \mathbb{N},
\emptyset, 2 \mathbb{N} + 1)$.

\begin{definition}\label{Models}
${\cal M}_1$ and ${\cal M}_2$ are of the form ${\cal M} = \langle
W, \leq, \mathbf{u}, D, \phi \rangle$; they are defined as
follows.
\begin{enumerate}
\item The base points for ${\cal M}_1$ and ${\cal M}_2$ are
$\mathbf{v}$ and $\mathbf{w}$ respectively; \item The set of
states for the models are as follows.
\begin{enumerate}
\item $W_1 = \{ (A,B,C) | \: \mathbf{v} \trianglelefteq (A,B,C),
\mbox{ where } B \cap \mathbf{v}_2 \mbox{ is infinite} \}$; \item
$W_2 = \{ (A,B,C) | \: \mathbf{w} \trianglelefteq (A,B,C), \mbox{
where } B \neq \emptyset \} \cup \{ \mathbf{w} \}$;
\end{enumerate}
\item The ordering on both ${\cal M}_1$ and ${\cal M}_2$ is
$\trianglelefteq$; \item $D_1 = D_2 = \mathbb{N}$; \item For $i =
1,2$, the values $\phi_i$ assigned to  the predicate symbols $P$
and $Q$ are defined by
\[ \phi_i(P) = \{ \langle v,a \rangle | a \in v_1 \cup v_2 \}, \]
\[ \phi_i(Q) = \{ \langle v,a \rangle | a \in v_1 \}. \]
\end{enumerate}
\end{definition}

The relation $\trianglelefteq$ is clearly reflexive and
transitive. Furthermore, if $(A,B,C) \trianglelefteq (D,E,F)$,
then because $F \subseteq C$, $A \cup B \subseteq D \cup E$. In
addition, $A \subseteq D$,  so that $\phi_i$ satisfies the
monotonicity property. Hence, both ${\cal M}_1$ and ${\cal M}_2$
are G-models.

\begin{lemma}\label{IJconditions}
\begin{enumerate}
\item The G-model ${\cal M}_1$ satisfies the condition $I(P,Q)$;
\item The condition $J(P,Q)$ fails in the G-model ${\cal M}_2$.
\end{enumerate}
\end{lemma}

{\bf Proof.} For the first part of the lemma, let $u = (A,B,C)$ be
a state in $W_1$, and $\mathbf{v}$ the base point of ${\cal M}_1$.
By definition, $B \cap \mathbf{v}_2$ is infinite. Choose $k \in  B
\cap \mathbf{v}_2$. Then $\mathbf{v} \Vdash_1 P\mathbf{k}$ and $u
\not\Vdash_1 Q\mathbf{k}$, showing that $I(P,Q)$ holds in ${\cal
M}_1$.

For the second part of the lemma, at the base point $\mathbf{w}$
of the G-model ${\cal M}_2$, we have  $\forall a \in D_2
(\mathbf{w} \Vdash_2 P{\bf a} \Rightarrow \mathbf{w} \Vdash_2 Q{\bf a})$,
because $\mathbf{w}_2 = \emptyset$. Hence, the condition $J(P,Q)$
fails in the model ${\cal M}_2$. \qed

\medskip

We next need to define a CD-asimulation $Z$ between the G-models
${\cal M}_1$ and ${\cal M}_2$. The following notations are useful
in stating the definition of $Z$ and in demonstrating its
properties. If $\vec{u}$ is a sequence of elements
from a set $D$, then $D \smallsetminus \vec{u}$ is the result of
removing all of the elements of $\vec{u}$ from $D$. For
$\vec{d},\vec{e} \in \mathbb{N}^{k}$, $S \subseteq \mathbb{N}$
we denote by $[\vec{d} \mapsto \vec{e}]$ the binary relation
$\{ \langle d_l,e_l\rangle\mid 1 \leq l \leq k \}$ and by
$[\vec{d} \mapsto \vec{e}]^{-1}(S)$ the set
$\{ d_l \:\: | \:\: e_l \in S, 1 \leq l \leq k \}$.

Because sequences $\vec{d},\vec{e} \in \mathbb{N}^{k}$ may contain
repetitions, the relation $[\vec{d} \mapsto \vec{e}]$ may not be a bijection,
and in fact, may not even be a function.
However, if $\vec{d}$ and $\vec{e}$ have the same pattern of repetitions
(for example, if $\vec{d} = 1,1,2,3,2$ and $\vec{e} = 4,4,5,6,5$) then
$[\vec{d} \mapsto \vec{e}]$ is a bijection, as is required in the next definition.

\begin{definition}\label{DefinitionofZ}
Relative to the models ${\cal M}_1$ and ${\cal M}_2$ given in
Definition \ref{Models}, the relation $Z$ is defined as follows:
\begin{enumerate}
\item $Z \subseteq \bigcup_{k \geq 0} [ (W_1 \times D_1^{k})
  \times (W_2 \times D_2^{k})] \cup [ (W_2 \times D_2^{k})
  \times (W_1 \times D_1^{k})]$;
\item $\langle (A,B,C), \vec{d} \rangle \: Z \: \langle (D,E,F), \vec{e} \rangle$,
  where $\vec{d} \in D_i^{k}$, $\vec{e} \in D_j^{k}$, $\{i,j\} = \{1,2\}$,
  if and only if the following conditions hold:
\begin{enumerate}
\item Relation $[\vec{d} \mapsto \vec{e}]$ is a bijection.

\item If $1 \leq l \leq k$ and  $d_l \in A$ then $e_l \in D$;

\item If $1 \leq l \leq k$ and $d_l \in B$, then $e_l \in D \cup
E$.
\end{enumerate}
\end{enumerate}
\end{definition}

It is easy to see that in the case when  both $\langle (A,B,C),
\vec{d} \rangle \: Z \: \langle (D,E,F), \vec{e} \rangle$ and
 $\langle (D,E,F), \vec{e} \rangle \: Z \: \langle (A,B,C), \vec{d}
\rangle$ hold, conditions 2(b),(c) of this definition are
equivalent, modulo other restrictions,  to the following ones:
\[
d_l\in A \text{ iff } e_l\in D;
\]
\[
d_l\in B \text{ iff } e_l\in E,
\]
for every $1 \leq l \leq k$.

\begin{lemma}\label{Zanasimulation}
 The relation $Z$ in Definition \ref{DefinitionofZ} is a CD-asimulation between the G-models ${\cal M}_1$ and ${\cal M}_2$.
\end{lemma}

{\bf Proof.} The first condition in Definition \ref{Asimulation}
is true by definition. For the second condition, assume that $v,
\vec{d} Z w, \vec{e}$ and that $v \forcesi P[\vec{\mathbf d}]$,
where $P(x_l)$ is atomic, $\vec{d} \in D_i^{k}$, and $1 \leq l
\leq k$. Thus we have $v \forcesi P[{\mathbf d}]$, where $d =
d_l$, so that $d \in v_1 \cup v_2$; it follows that $e = e_l \in
w_1 \cup w_2$, by Definition \ref{DefinitionofZ}, showing that $v
\forcesj P[\vec{\mathbf e}]$. The proof for atomic formulas
$Q(x_l)$ is similar.

For the third condition in Definition \ref{Asimulation}, assume
that $t, \vec{d} Z u, \vec{e}$, where $t = (A,B,C)$, $u =
(D,E,F)$, and $u \leq_j v$, $v = (G,H,I)$. By definition, $u
\trianglelefteq v$. Two cases arise here: $B$ is infinite, or $B =
\emptyset$.

In the first case, define $w = (J,K,L)$ as follows:
\begin{eqnarray*}
 J & = & (A \smallsetminus \vec{d}) \cup [\vec{d} \mapsto \vec{e}]^{-1}(G) ;\\
 K & = & (B \smallsetminus \vec{d}) \cup [\vec{d} \mapsto
 \vec{e}]^{-1}(H);\\
 L & = & (C \smallsetminus \vec{d}) \cup [\vec{d} \mapsto \vec{e}]^{-1}(I).
\end{eqnarray*}
We need to show that $(J,K,L)$ is indeed a
successor of $(A,B,C)$ in ${\cal M}_i$, i.e. that $(A,B,C)
\trianglelefteq (J,K,L)$ and that if  $B \cap \mathbf{v}_2$ is
infinite, then so is $K \cap \mathbf{v}_2$. For the latter claim
note that if $B \cap \mathbf{v}_2$ is infinite, then so is $(B
\smallsetminus \vec{d}) \cap \mathbf{v}_2$, given that $\vec{d}$
is finite. But since, according to our definition of $K$ we have
$(B \smallsetminus \vec{d}) \subseteq K$, $K \cap \mathbf{v}_2$
must be infinite, too, showing that $(J,K,L)$ is a quasi-partition
in ${\cal M}_i$.

The former claim, that is to say, the claim that $A \subseteq J$
and $L \subseteq C$, can be verified as follows. If $a \in A$, and
$a$ is in $\vec{d}$, say $a = d_l$, then $e_l \in D$, so $e_l \in
G$, from which it follows that $d_l = a \in J$. Since $I \subseteq
F$, we have $[\vec{d} \mapsto \vec{e}]^{-1}(I) \subseteq [\vec{d}
\mapsto \vec{e}]^{-1}(F) \subseteq C$, showing that $L \subseteq
C$. Finally, it follows by construction that for $1 \leq l \leq
k$, $d_l \in J \Leftrightarrow e_l \in G$, and $d_l \in K
\Leftrightarrow e_l \in H$. Hence, $w, \vec{d} Z v, \vec{e}$ and
$v, \vec{e} Z w, \vec{d}$, completing the proof of the third
condition in the first case.

In the second case, where $B = \emptyset$, we have $(A,B,C) =
\mathbf{w} = ( 2 \mathbb{N}, \emptyset, 2 \mathbb{N} + 1)$. In
this case, we partition $C \smallsetminus \vec{d}$ into two
disjoint infinite sets $C_1$ and $C_2$, and define $w = (J,K,L)$
as follows:
\begin{eqnarray*}
 J & = & (A \smallsetminus \vec{d}) \cup [\vec{d} \mapsto \vec{e}]^{-1}(G) ;\\
 K & = & C_1 \cup [\vec{d} \mapsto \vec{e}]^{-1}(H);\\
 L & = & C_2 \cup [\vec{d} \mapsto \vec{e}]^{-1}(I).
\end{eqnarray*}
Then $J,K,L$ are all infinite by construction. The remainder of
the proof of the third condition in the second case is essentially
the same as for the first case, the only difference being that in
this case we do not need to establish that if  $B \cap
\mathbf{v}_2$ is infinite, then so is $K \cap \mathbf{v}_2$.

For the fourth condition in Definition \ref{Asimulation}, assume
that $t \in W_i$, $(t,\vec{d} Z u, \vec{e})$, where $t = (G,H,I)$,
$u = (J,K,L)$, and $f \in D_i$. If $f = d_l$ for some $1 \leq l
\leq k$, then set $g := e_l$. Otherwise, given that $f \notin
\vec{d}$, choose $g \in D_j$ so that $g \in J \smallsetminus
\vec{e}$ (this is possible because $J$ is infinite). Thus we have
shown the existence of a $g$ so that $t,\vec{d},f Z u, \vec{e},g$.

The fifth condition is proved by a similar argument. Namely,
assume that $t \in W_i$, $(t,\vec{d} Z u, \vec{e})$, where $t =
(G,H,I)$, $u = (J,K,L)$, and $g \in D_j$. If $g = e_l$ for some $1
\leq l \leq k$, then set $f := d_l$. Otherwise, since $g \notin
\vec{e}$, choose $f \in D_i$ so that $g \in  I \smallsetminus
\vec{d}$ (this is possible because $I$ is infinite).  \qed

\medskip

We now have all the material required to refute interpolation for \cd.

\begin{theorem}
The implication $\Gamma \rightarrow \Delta$ is valid in all \G-models,
but has no interpolant in \cd.
\end{theorem}

{\bf Proof.} The validity of the implication is proved in Lemma
\ref{Implicationlemma}. Now let us assume that the implication has
an interpolant in \cd, that is to say, a sentence $\Theta$ in the
language $L(P,Q)$ so that both $\Gamma \rightarrow \Theta$ and
$\Theta \rightarrow \Delta$ are valid in all G-models. Then
$\Theta$ is also an interpolant for the second-order implication
$\exists R \Gamma \rightarrow \forall S \Delta$.

Since the G-model ${\cal M}_1$ satisfies the condition $I(P,Q)$,
by Lemma \ref{IJconditions}, its base point also satisfies the
second-order condition $\exists  R \Gamma$, by Lemma
\ref{Semanticalinterpolants}, so $\mathbf{v} \Vdash_1 \Theta$.
Moreover, it is clear from Definition \ref{DefinitionofZ} that we
have $\langle \mathbf{v}, \Lambda \rangle \: Z \: \langle
\mathbf{w}, \Lambda \rangle$, where $\Lambda$ is the empty
sequence. Lemma \ref{Zanasimulation} shows that $Z$ is a
CD-asimulation between ${\cal M}_1$ and ${\cal M}_2$, so by Lemma
\ref{Soundness}, $\mathbf{w} \Vdash_2 \Theta$, hence $\mathbf{w}
\Vdash_2 \forall S\Delta$. However, Lemma
\ref{Semanticalinterpolants} shows that ${\cal M}_2$ must satisfy
$J(P,Q)$, contradicting Lemma \ref{IJconditions}. Consequently, no
such interpolant for the implication $\Gamma \rightarrow \Delta$
can exist.
\qed

\bibliographystyle{habbrv}

\bibliography{logic}

\begin{thebibliography}{10}

\bibitem{Beth1956}
E.~W. Beth.
\newblock Semantic construction of intuitionistic logic.
\newblock {\em Koninklijke Nederlandse Akademie van Wetenschappen,
  Mededelingen, Nieuwe Reeks}, 19:357--388, 1956.

\bibitem{Beth1959}
E.~W. Beth.
\newblock {\em The Foundations of Mathematics: a study in the philosophy of
  science}.
\newblock North Holland Publishing Company, Amsterdam, 1959.

\bibitem{BlackburndeRijkeVenema2001}
P.~Blackburn, M.~de~Rijke, and Y.~Venema.
\newblock {\em Modal Logic}.
\newblock Cambridge University Press, 2001.
\newblock Cambridge Tracts in Theoretical Computer Science 53.

\bibitem{Fine1979}
K.~Fine.
\newblock Failures of the interpolation lemma in quantified modal logics.
\newblock {\em Journal of Symbolic Logic}, 44:201--206, 1979.

\bibitem{Gabbay1969}
D.~Gabbay.
\newblock {Montague type semantics for nonclassical logics. I}.
\newblock Technical report, Hebrew University of Jerusalem, 1969.

\bibitem{Gabbay1977}
D.~Gabbay.
\newblock Craig interpolation theorem for intuitionistic logic and extensions
  {Part III}.
\newblock {\em Journal of Symbolic Logic}, 42:269--271, 1977.

\bibitem{GabbayMaksimova}
D.~Gabbay and L.~Maksimova.
\newblock {\em Interpolation and Definability: Modal and Intuitionistic
  Logics}.
\newblock Oxford University Press, 2005.

\bibitem{Gornemann1971}
S.~G{\"o}rnemann.
\newblock A logic stronger than intuitionism.
\newblock {\em Journal of Symbolic Logic}, 36:249--261, 1971.

\bibitem{Grzegorczyk1964}
A.~Grzegorczyk.
\newblock A philosophically plausible formal interpretation of intuitionistic
  logic.
\newblock {\em Indagationes Mathematicae}, 26:596--601, 1964.

\bibitem{Klemke1971}
D.~Klemke.
\newblock {Ein Henkin-Beweis f{\"u}r die Vollst{\"a}ndigkeit eines Kalk{\"u}ls
  relativ zur Grzegorczyk-Semantik}.
\newblock {\em Archiv f{\"u}r Mathematische Logik und Grundlagenforschung},
  14:148--161, 1971.

\bibitem{Kripke1965}
S.~A. Kripke.
\newblock {Semantical analysis of intuitionistic logic. I}.
\newblock In {\em Formal Systems and Recursive Functions}, pages 92--130, 1965.
\newblock Proceedings of the Eighth Logic Colloquium, Oxford, July 1963.

\bibitem{Kripke1983}
S.~A. Kripke.
\newblock Review of \cite{Fine1979}.
\newblock {\em Journal of Symbolic Logic}, 44:486--488, 1983.

\bibitem{LopezEscobar1981}
E.~L\'{o}pez-Escobar.
\newblock On the interpolation theorem for the logic of constant domains.
\newblock {\em Journal of Symbolic Logic}, 46:87--88, 1981.

\bibitem{LopezEscobar1983}
E.~L\'{o}pez-Escobar.
\newblock {A second paper ``On the interpolation theorem for the logic of
  constant domains''}.
\newblock {\em Journal of Symbolic Logic}, 48:595--599, 1983.

\bibitem{Mints2000}
G.~Mints.
\newblock {\em A Short Introduction to Intuitionistic Logic}.
\newblock Kluwer Publishers, 2000.

\bibitem{Olkhovikov2012}
G.~K. {Olkhovikov}.
\newblock {Model-theoretic characterization of predicate intuitionistic
  formulas}.
\newblock {\em ArXiv e-prints}, Feb. 2012, 1202.1195.

\bibitem{Takeuti1987}
G.~Takeuti.
\newblock {\em Proof Theory}.
\newblock North-Holland, 1987.

\bibitem{VanBenthem2010}
J.~van Benthem.
\newblock Frame correspondences in modal predicate logic.
\newblock In {\em Proofs, Categories and Computations: Essays in Honor of
  Grigori Mints}. College Publications, London, 2010.

\end{thebibliography}

\end{document}